\documentstyle[12pt]{article}

\def\1{\'{\i}}

\def\be{\begin{equation}}
\def\ee{\end{equation}}
\def\bea{\begin{eqnarray}}
\def\eea{\end{eqnarray}}

\def\Osc{H_4}  

\def\R{{\rm I\kern-.2em R}}

\def\aa{N} 
\def\ap{{A_+}} 
\def\am{{A_-}}
\def\bb{M}

\def\aaa{n} 
\def\aap{{a_+}}
\def\aam{{a_-}}
\def\bbb{m}

\def\haaa{{\hat n}} 
\def\haap{{\hat{a}_+}}
\def\haam{{\hat{a}_-}}
\def\hbbb{\hat{m}}

 \def\1{\'{\i}}

\begin{document}

\ 
\vspace{2cm}

\begin{center}
 {\LARGE{\bf{Harmonic Oscillator Lie Bialgebras}}} 
  
 {\LARGE{\bf{and their Quantization}}} 
\end{center}

\bigskip

\begin{center} Angel Ballesteros and  Francisco J. Herranz
\end{center}

\begin{center} {\it {  Departamento de F\1sica, Universidad
de Burgos} \\   Pza. Misael Ba\~nuelos, \\
E-09001, Burgos, Spain}
\end{center}

\begin{abstract}
All possible Lie bialgebra structures on the harmonic oscillator
algebra are explicitly derived and it is shown that all of them are of the
coboundary type. A non-standard quantum oscillator is introduced as a
quantization of a triangular Lie bialgebra, and  a universal $R$-matrix
linked to this new quantum algebra is presented.
\end{abstract}

\section{Introduction}

The theory of Hopf algebra deformations of   Universal Enveloping Algebras
is, by construction, intimately linked to Lie bialgebras (a detailed
exposition can be found in \cite{PC}). In fact, once a uniparametric
deformation of $Ug$ is given, a unique Lie bialgebra structure $(g,\delta)$
on $g$ is defined via
\be
\delta:=(\Delta_z - \sigma\circ \Delta_z)\qquad \mbox{mod}\ z^2,
\label{zz}
\ee
where $\sigma$ is the permutation map   $\sigma\,(X\otimes Y)=Y\otimes X$,
$\Delta_z$ is the deformed coproduct and  $z$ the deformation parameter.  As
it was shown by Drinfel'd
\cite{Dr}, the cocommutator $\delta$  characterizes completely the
Poisson--Lie group whose quantization is given by the Hopf algebra dual to
$U_z g$. Therefore, Lie bialgebra structures on $g$ are the natural
candidates to explore the zoo of quantum algebra deformations of a given
algebra.  Note that (\ref{zz}) means that we have extracted the first order
term in $z$ by preserving such a deformation parameter as a multiplicative
factor within
$\delta$.  The reason for this is to open the possibility of considering Lie
bialgebras corresponding to multiparameter deformations of $U g$, that will
contain simultaneously various deformation parameters. 
 In this context two natural questions arise: 

\noindent $1.$ Given a precise Lie algebra, is it possible to
obtain explicitly all its Lie bialgebra structures?

\noindent $2.$ Once $\delta$
is given, can its quantum deformation $U_z g$ be constructed?

Essentially, solving both problems will be equivalent to
the obtention of all quantum algebra deformations
of $g$. This paper is mainly devoted to answer in the affirmative both
questions for the physically meaningful  case of the harmonic oscillator
algebra
$h_4$.   

In particular, in Section 2 we explicitly  compute all Lie bialgebra
structures on $h_4$, the main result of the paper being the proof of the
coboundary character of all of them. This fact   means that the
classification of the coboundary Lie bialgebras of $h_4$ given in \cite{osc}
provides a full characterization of all posible quantum  (multiparametric)
deformations of this algebra. Moreover, all their corresponding coproducts
were constructed in
\cite{osc} by using a link between Lie  bialgebras and a general method for
the construction of coassociative coproducts on an arbitrary associative
algebra \cite{Lyak} (recall that a systematic  derivation of quantum
oscillator groups was given in \cite{Vero}). Finally, in Section 3 we shall
concentrate our attention onto a precise non-standard quantization of the
oscillator algebra that can be properly called the ``Jordanian
$q$-oscillator". Its Hopf algebra structure  will be analysed, and the
universal quantum $R$-matrix presented for this triangular quantization.


\section{Harmonic oscillator Lie bialgebras}

We begin by recalling that a Lie bialgebra $(g,\delta)$ is a Lie algebra $g$
endowed with a map $\delta:g\to g\otimes g$ such that

\noindent i) $\delta$ is a 1-cocycle, i.e.,
\be
\delta([x,y])=[\delta(x),\, 1\otimes y+ y\otimes 1] + 
[1\otimes x+ x\otimes 1,\, \delta(y)], \quad \forall x,y\in
g. 
\label{bb}
\ee
\noindent ii) The dual map $\delta^\ast:g^\ast\otimes g^\ast \to
g^\ast$ is a Lie bracket on $g^\ast$.

A Lie bialgebra $(g,\delta)$ is called a coboundary Lie bialgebra if there
exists an element $r\in g\otimes g$ (the classical $r$-matrix), such that 
\be
\delta(x)=[1\otimes x + x \otimes 1,\,  r], \qquad  \forall x\in
g.
\label{bc}
\ee
When the $r$-matrix is such that its Schouten bracket vanishes, we shall
say that $(g,\delta (r))$ is a {\em non-standard} (or triangular) Lie
bialgebra. On the contrary, we shall have a {\em standard}
structure.
Finally, two Lie bialgebras $(g,\delta)$  and $(g,\delta')$ are said to be
equivalent if there exists an automorphism $O$ of $g$ such that
$\delta'=(O\otimes O)\circ\delta\circ O^{-1}$. 

Let us now consider the specific case of $h_4$, with generators
$\{\aa,\ap,\am,\bb\}$ and commutation rules \be
[\aa,\ap]=\ap,\quad [\aa,\am]=-\am,\quad [\am,\ap]=\bb,\quad 
[\bb,\cdot\,]=0 .
\label{aa}
\ee
The most general cocommutator  $\delta:h_4\rightarrow h_4 \otimes h_4$ will
be a linear combination (with real coefficients) \be
\delta(X_i)=f_i^{jk}\,X_j\wedge X_k,
\label{gen}
\ee
of skew-symmetric products of  the generators $X_l$ of the algebra. Such a
completely general cocommutator has to be computed by firstly imposing the
cocycle condition (\ref{bb}) onto an arbitrary expression (\ref{gen}). This
leads to the following six-parameter (pre)cocommutator:
\bea && \delta(N)= a_1\,N\wedge A_+ 
+ a_2\,N\wedge A_-  +  a_5\,A_+\wedge M  + a_6\,A_-\wedge M,\cr
&& \delta(A_+)=
a_2\,N\wedge M  + a_2\,A_+\wedge A_-  +
a_3\,A_+\wedge M,\cr
&& \delta(A_-)=
a_1\,N\wedge M  - a_1\,A_+\wedge A_-  +
a_4\,A_-\wedge M,\cr
&& \delta(M)=0.
\label{is}
\eea 
Afterwards, Jacobi identity has to be imposed onto $\delta^\ast: h_4^\ast
\otimes h_4^\ast \rightarrow h_4^\ast $  in order to guarantee that this map
defines a Lie bracket. It is easy to  check that this condition leads to the
set of equations \be a_1\,a_3=0,\qquad  a_1\,a_2=0,\qquad 
a_2\,a_4=0.
\label{dual}
\ee
The set of solutions can be splitted  into three disjoint classes and
they give rise to the following complete classification of Lie bialgebra
structures of
$h_4$:

\noindent  Type A bialgebras: $a_1\neq 0$, $a_2=0$ and 
$a_3=0$. 

\noindent  Type B bialgebras: $a_1= 0$, $a_2\neq 0$ and 
$a_4=0$. 

\noindent  Type C bialgebras: $a_1= 0$ and $a_2=0$. 

Note that, from (\ref{is}), all  these three types are four-parameter
families of Lie bialgebras in which the  cocommutator of the central
generator $M$ always vanishes.  Moreover for types A and B, respectively,
$\delta(A_+)=0$ and $\delta(A_-)=0$.

As we have mentioned in the  introduction, the classification of all
coboundary Lie bialgebra structures of $h_4$ was developed in \cite{osc}.
With this in mind, it is inmmediate to show that the classification just
obtained coincides exactly with that of the coboundary structures, i.e.,
non-coboundary oscillator Lie bialgebras do not exist. In fact, Type A
bialgebras coincide with Type I$_+$ ones in \cite{osc} under the following
identification of the parameters
\be
a_1=\alpha_+\neq 0,\qquad a_4=2\, \vartheta,\qquad a_5=\beta_+,\qquad
a_6=-\beta_- .
\label{iduno}
\ee
The additional condition
\be
4\,a_1\,a_6 + a_4^2=0,
\label{iddos}
\ee
gives rise to non-standard (or triangular) Lie bialgebras. The classical
$r$-matrix for Type A structures is 
\be
r_A=a_1\,N\wedge A_+ + a_4\,(N\wedge M - A_+\wedge A_-)/2 + 
a_5\, A_+\wedge M - a_6\, A_-\wedge M. 
\label{rmatA}
\ee

The connection between Type B and Type I$_-$ structures is also
straightforward if we put
\be
a_2=-\alpha_-\neq 0,\qquad a_3=-2\,\vartheta, \qquad a_5=\beta_+,\qquad
a_6=-\beta_-,
\label{idtres}
\ee
and non-standard Lie bialgebras of Type B appear when
\be
4\,a_2\,a_5 + a_3^2=0.
\label{idcuatro}
\ee
The most general classical $r$-matrix for this kind of Lie bialgebras is 
\be
r_B= - a_2\,N\wedge A_- - a_3\,(N\wedge M + A_+\wedge A_-)/2 + 
a_5\, A_+\wedge M - a_6\, A_-\wedge M. 
\label{rmatB}
\ee

This Type B family is easily proven to be equivalent (as a multiparametric
Lie bialgebra) to the Type A one by means of the automorphism of $h_4$ that
interchanges the $A_+$ and $A_-$ generators. Finally, if we write Type C
bialgebras by using that
\be
a_3=-\vartheta - \xi,\qquad a_4=\vartheta - \xi,\qquad a_5=\beta_+,\qquad
a_6=-\beta_-, \label{idcinco}
\ee
the result is just the cocommutator of the Type II bialgebras in \cite{osc}.
Moreover, in case that
\be
a_3+a_4=0,
\label{idseis}
\ee
we shall have recovered the non-standard cases. The classical  Type C
$r$-matrix is 
\be
r_C= (a_4 - a_3) \, N\wedge M/2 - (a_4 + a_3) \, A_+\wedge A_-/2 +
a_5\, A_+\wedge M - a_6\, A_-\wedge M. 
\label{rmatC}
\ee
 
The obtention of a coassociative quantum  coproduct quantizing each of these
Lie bialgebras was developed in \cite{osc}.  In some cases, automorphisms of
the oscillator algebra are needed to simplify the Lie bialgebras before
quantizing them. For instance, by using the automorphism
\be
\aa'=\aa- (a_5/a_1)\bb,
\label{ggi}
\ee
of $h_4$ it can be shown how $a_5$ is
an irrelevant parameter in Type A bialgebras. Likewise   $a_6$ plays 
a trivial role in the family B.


\section{The Jordanian quantum oscillator}

The (standard) quantum oscillator algebra obtained in
\cite{Celeghinidos,Sierra} by a contraction method can  be easily recovered
from our classification. In fact, the classical $r$-matrix underlying that
deformation is 
\cite{Celeghinidos}:
\be
r= - z\,(\aa\otimes\bb + \bb\otimes\aa) + 2z \,\am\otimes \ap,
\label{pa}
\ee
and its skew-symmetric ($r_-$) part
\be
r_-=(r-\sigma\circ r)/2=
z \, \am\wedge \ap   
\label{pc}
\ee
is, in our notation, a
standard classical $r$-matrix of Type C with parameters
$a_5=a_6=0$ and $a_3=a_4=z$. The associated oscillator
bialgebra reads:
\be
\delta(\aa)=\delta(\bb)=0,\qquad \delta(\ap)=z\, \ap\wedge\bb,\qquad
\delta(\am)=z\, \am\wedge\bb .
\label{pd}
\ee

At this point, it is worth recalling that this  deformation was firstly
obtained in \cite{Celeghinidos} by contracting the standard deformation of
$sl(2,\R)$. Therefore, it seems pertinent to wonder whether a non-standard
deformation of the oscillator algebra related to the so-called  Jordanian
quantum $sl(2,\R)$ \cite{Demidov,Zakr,Ohn} exists. This question  can be
easily addressed at the Lie bialgebra level with the help of the
clasification here presented as follows.

The classical $r$-matrix underlying the  Jordanian deformation of $sl(2,\R)$
is
\be
r=z\, J_3\wedge J_+ .
\label{jor}
\ee
Now, if we examine the non-standard oscillator Lie bialgebra of Type A
with $a_4=a_5=a_6=0$ and $a_1= z$, its classical $r$-matrix would be
\be 
 r=z\, \aa\wedge \ap .
\label{ea}
\ee 
By taking into account that $\{N,A_+\}$  define the same Lie algebra as
$\{J_3,J_+\}$, the deformation of this positive Borel subalgebra induced by
both (\ref{jor}) and (\ref{ea}) will be the same. Moreover, the full
non-standard deformation of $sl(2,\R)$ can be obtained from the information
contained in (\ref{jor}). Now it is immediate to compute the harmonic
oscillator cocommutator derived from (\ref{ea}):
\bea
&&\delta(\ap)=0,\quad \delta(\bb)=0,\quad 
  \delta(\aa)=z\, \aa\wedge\ap  ,\cr
&& \delta(\am)=z\, (\am\wedge \ap+\aa\wedge\bb )   .
\label{eb}
\eea
The quantum algebra obtained by quantizing  this Lie bialgebra can therefore
be properly called the ``Jordanian $q$-oscillator", and has the following
coproduct and commutation rules
\bea
&&\Delta(\ap)=1\otimes \ap +\ap \otimes 1,\cr
&& \Delta(\bb)=1\otimes \bb +\bb \otimes 1,\cr
&&\Delta(\aa)=1\otimes \aa +\aa \otimes e^{z\ap},\cr 
&&\Delta(\am)=1\otimes \am +\am \otimes e^{z\ap}+z\aa\otimes \bb e^{z\ap} ;
\label{ed}
\eea
\be 
 [\aa,\ap]=\frac{e^{z\ap}-1}{z},\quad [\aa,\am]=-\am
,\quad [\am,\ap]=\bb e^{z\ap},\quad 
[\bb,\cdot\,]=0 .
\label{eg}
\ee 

An essential feature of this deformation is the existence of a quantum
universal $R$-matrix, that can be derived by taking into account again that 
$\aa$ and $\ap$ generate a deformed Hopf subalgebra. Such a quantum Hopf
subalgebra does have a universal $R$-matrix, already derived in 
\cite{Tmatrix}:
\be
R=\exp\{-z\ap\otimes\aa\}\exp\{z\aa\otimes\ap\} .
\label{ei}
\ee
moreover, the solution (\ref{ei}) of the quantum  Yang--Baxter equation can
be shown \cite{osc} to fulfill the relations
\be
\sigma\circ \Delta(X)=R\Delta(X)R^{-1},\qquad
\mbox{for}\ \ X\in\{\aa,\ap,\am,\bb\} .
\label{ej}
\ee
Therefore, (\ref{ei}) is also a universal $R$-matrix for the Jordanian
$q$-oscillator. 

As a consequence, the construction of the corresponding 
quantum oscillator group is thus available by following the FRT approach
\cite{FRT}. For that purpose we need to recall  that a matrix representation
of the classical oscillator algebra  is provided by  
\bea
&&D(\aa)=\left(\begin{array}{ccc}
 0 &0 & 0 \\ 0 & 1 & 0 \\ 0 & 0 & 0 
\end{array}\right),\qquad  
D(\ap)=\left(\begin{array}{ccc}
 0 &0 & 0 \\ 0 & 0 & 1 \\ 0 & 0 & 0 
\end{array}\right),\cr 
&&D(\am)=\left(\begin{array}{ccc}
 0 &1 & 0 \\ 0 & 0 & 0 \\ 0 & 0 & 0 
\end{array}\right),\qquad 
D(\bb)=\left(\begin{array}{ccc}
 0 &0 & 1 \\ 0 & 0 & 0 \\ 0 & 0 & 0 
\end{array}\right).
\label{ac}
\eea
A generic element $T^D$ of the oscillator group   $\Osc$
coming from this representation is: 
\bea
&&T^D=\exp\{\bbb D(\bb)\}\exp\{\aam D(\am)\}
\exp\{\aap D(\ap)\}\exp\{\aaa D(\aa)\}\cr
&&\quad =\left(\begin{array}{ccc}
 1 &\aam e^\aaa & \bbb+\aam\aap \\ 0 & e^\aaa & \aap \\ 0 & 0 & 1 
\end{array}\right) .
\label{ad}
\eea
By taking into account that the representation (\ref{ac})   is also a
realization of the deformed algebra (\ref{eg}), the universal $R$-matrix
(\ref{ei}) is represented as
\be
D(R)=I\otimes I + z(D(\aa)\otimes D(\ap) - D(\ap)\otimes D(\aa)), 
\label{rr}
\ee
($I$ is the $3\times 3$
identity matrix). Finally, by considering  non-commutative entries in
(\ref{ad}), the following quantum oscillator group is obtained (quantum
coordinates are denoted by hats):   
\bea
&&\Delta( \haaa)=1\otimes  \haaa +  \haaa \otimes 1,\cr
&&\Delta( \haap)=e^{ \haaa}\otimes  \haap + \haap \otimes 1,\cr
&&\Delta( \haam)=e^{- \haaa}\otimes  \haam + \haam \otimes 1,\cr
&&\Delta( \hbbb)=1\otimes  \hbbb + \hbbb \otimes 1
-e^{- \haaa} \haap\otimes  \haam ;
\label{ek}
\eea
\bea
&&[ \haaa, \haap]=z\, (e^{  \haaa} -1),\qquad
[ \haaa, \haam]=0,\qquad
[ \haam, \haap]=z\,  \haam,\cr
&&[ \haaa, \hbbb]=z\,  \haam,\qquad
[ \haap, \hbbb]=z\, \haam \haap ,\qquad
 [ \haam, \hbbb]=-z\,{ \haam}^2  .
\label{en}
\eea
This Jordanian quantum oscillator group is just  a Weyl quantization of  the
Poisson--Lie bracket (expressed in local  coordinates) on the classical
oscillator group defined by the $r$-matrix  (\ref{ea}) via the Sklyanin
bracket. 

We would like to point out some interesting  features of this new Hopf
algebra deformation of the oscillator algebra.  Among them we emphasize
that, in contradistinction to the standard deformation,  exponentials of
$A_+$ are the essential constituents of the deformed Hopf algebra. In the
standard deformation, this role was played by the mass $M$, that now is also
primitive and central. Another central element, the quantum Casimir, can be
computed and reads
\be
C_z=2\aa\bb+  
\frac{e^{-z\ap}-1}{z}\, \am +\am \, \frac{e^{-z\ap}-1}{z}.
\label{eh}
\ee
The eigenvalues of this Casimir can be used to label  the representation
theory of this non-stan\-dard $q$-oscillator, that presents stringent
differences with respect to that of the standard deformation and also in
relation to the well known $q$-oscillator of Biedenharn and Macfarlane
\cite{Bi,Mf} (we recall that the latter is not a Hopf algebra deformation,
therefore it cannot be recovered from our classification). For instace, the
Jordanian $q$-oscillator can be realized  in terms of ordinary boson
operators
$a_-$, $a_+$, verifying  $[a_-,a_+]=1$, in the form
\bea
&& \ap=\aap,\qquad  \bb=\delta 1 ,\cr
&& \am=\delta e^{z\aap}\aam +  \delta\beta\frac{z}2\, 
e^{z\aap},\label{gtq}\\ 
&&\aa= \frac{e^{z\aap}-1}{z}\,\aam+\beta \,
\frac{e^{z\aap}+1}{2},\nonumber
\eea
where $\beta$, $\delta$ are related to the eigenvalue of the Casimir as
$C_z=\delta(2\beta -1)$.   A systematic study of
the representation theory of non-standard deformations,  as well as a
derivation of the Hopf algebra structure of the Jordanian $q$-oscillator
starting from a non-standard quantum deformation of a pseudo-extended
$sl(2,\R)$ algebra can be found in \cite{BHN}.


\section*{Acknowledgments}
This work has been
partially supported by DGICYT (Project PB94--1115) from the
Ministerio de Educaci\'on y Ciencia de Espa\~na.

\end{document}